\documentclass[12pt]{article}
\usepackage{color}

\def\noi{\noindent}

\def\skip{\vspace{2mm}}
\newcount\u
\def\Remark{\skip\noi{{\bf{Remark \number\u.}}} \advance\u by 1}

\usepackage{amsmath,amssymb,amscd}

\usepackage{graphics}                   
\usepackage[dvips]{epsfig}
                   
\usepackage{delarray,hhline,latexsym}

\usepackage{float}

\usepackage{enumerate}

\usepackage{fancyhdr}

\def\noi{\noindent}

\newtheorem{props}{Proposition}

\newcounter{defnctr}
\newtheorem{defns}[defnctr]{Definition}

\newtheorem{corrs}[props]{Corollary}

\newtheorem{lemms}[props]{Lemma}

\newcounter{exctr}
\newtheorem{es}[exctr]{Example}

\newtheorem{ess}[exctr]{Examples}

\begin{document}

\title{Paradoxical decompositions and finitary rules}
\author {Robert Samuel Simon and Grzegorz Tomkowicz}

\maketitle
\thispagestyle{empty}
\vfill

\noi  Robert Simon 
\newline London School of Economics
\newline 
Department of Mathematics\newline 
Houghton Street\newline 
 London WC2A 2AE\newline                    
email: R.S.Simon@lse.ac.uk
\vskip.5cm

 \noi Grzegorz Tomkowicz\\
  Centrum Edukacji $G^2$\\
        ul.Moniuszki 9\\
        41-902 Bytom\\
        Poland\\
      e-mail: gtomko@vp.pl

\vfill

\date{}

\setcounter{page}{-1}
  \noi

\newpage
\vskip2cm
\thispagestyle{empty}
\begin{center}
{\bf Abstract}
\end{center}
\vskip1cm

\noi  We  colour every point $x$  of a probability space $X$
according   to  the
colours of a finite list  $x_1, x_2, \dots , x_k$ of points
 such that each of the  $x_i$,  as a function of $x$, is 
 a measure preserving transformation.  We ask two questions
 about  a colouring  rule: 
 (1) does there exist a finitely additive extension
 of the probability measure for which the $x_i$
 remain measure preserving and also a  colouring obeying the rule almost
  everywhere that is 
  measurable with respect to this extension?,
 and (2) does there exist any colouring obeying the rule almost everywhere?
 If the answer to the first question is no and to the second question yes, 
 we say that the colouring  rule is paradoxical. A paradoxical colouring 
 rule not only allows for a paradoxical partition of the space,
 it  requires one.      
 We pay special attention to generalizations of the Hausdorff paradox.

\vskip2cm

\noi {\bf Key words}: measurable  colouring, Cayley graphs, 
Hausdorff Paradox, paradoxes in measure theory

   \newpage

\section{Introduction} 
H. Lebesgue asked if his measure can be defined on all subsets of the real line $\mathbb{R}$. A negative answer was provided by G. Vitali.
Lebesgue measure $\lambda$ on $\mathbb{R}$ is a countably additive measure that is isometry invariant
 and such that $\lambda([0,1]) = 1$. Therefore, in the light of Vitali's construction, it is natural to ask if one can relax the condition of
 countable additivity and demand that such a measure to be defined on all subsets of $\mathbb{R}$. F. Hausdorff showed in [3]
 that there is no finitely additive rotation-invariant measure $\mu$ defined on all subsets of the unit sphere $\mathbb{S}^2$
 and such that $\mu(\mathbb{S}^2) = 1$. His construction can be easily modified to obtain that there is no isometry invariant finitely additive measure
 $\mu$ defined on all subsets of $\mathbb{R}^n$ ($n \geq 3$) and such that $\mu([0,1]^n) = 1$ (see [6] for a simple proof). Finally, S. Banach [1] proved that there is a finitely additive measure with required properties in $\mathbb{R}^n$ for $n=1,2$ (see [6], Cor. 12.9).  Both Vitali's and Hausdorff's constructions use Axiom of Choice (AC) for uncountable families of sets.
 Namely, the authors construct some sets which then 
 satisfy some rules (related to  partitions of the space) that violate conditions required by the measures in question.\\
 \indent In the present paper
 we reverse the above  perspective. We consider how to partition a probability space $X$ into
  subsets $A_1, \dots, A_k$ according to rules
 defined with   a  group or  semi-group  $G$ acting on the 
 space  $X$ in a  measure preserving way.
 We call such rules {\em colouring} rules, and we are interested  in colouring
  rules
 where  the \emph{only}  sets that satisfy these rules
 cannot be measurable with respect to \emph{any} finitely additive
 $G$-invariant measure $\mu$ extending the initial
  probability measure such that $\mu(X) = 1$.
  Furthermore we are not really
   interested in  such   rules that   are merely inconsistent,
  but rather in  those for
  which, using AC, the rule    can
  be satisfied in {\emph some} way. We call colouring
  rules  {\em paradoxical} when they can be satisfied
   almost everywhere but never in the finitely additive
   measurable  way described above.  \\  
   \indent Let us analyze closer the Hausdorff's construction in
           [3] that the unit sphere $\mathbb{S}^2$ can be decomposed,
           modulo a countable set $D$ (the intersection with
            the sphere of the
            set of all axes
            from the  appropriate subgroup
            of rotations)
            into the sets $A, B$ and $C$ satisfying the relation\\
 
 $(*)$ $A \cong B \cong C \cong B \cup C$,\\ 

  \noindent where $\cong$ denotes the congruence of sets. Here the congruences are witnessed by a subgroup of   $SO_3(\mathbb{R})$ isomorphic to $\mathbb{Z}_2 * \mathbb{Z}_3$ (and for
   convenience we will call  this subgroup  $\mathbb{Z}_2 * \mathbb{Z}_3$).  With $\sigma$ generating $\mathbb{Z}_2$
   and $\tau$ generating $\mathbb{Z}_3$, we have the set
   of rules $\sigma (A) = B\cup C$,
   $\sigma (B\cup C)=A$, $\tau (A) = B$, $\tau (B)=C$, and $\tau (C) =A$. We call these rules
   the {\em Hausdorff colouring rules}.
   Measurability with respect to any measure (where  both $\sigma$ and $\tau$ are still 
    measure preserving)  would  require
   that both   half and one-third  of the space must be taken by the set $A$. And  from these three sets, we can also   create six
   sets that partition the whole
    space, namely
    $A\cap \sigma (B)$, $A\cap \sigma (C)$, $\tau (A\cap \sigma (B))$,  $\tau (A\cap \sigma (C))$,  $\tau^2 (A\cap \sigma (B))$ and
    $\tau^2 (A\cap \sigma (C))$, such that by  choosing the  appropriate
     rotations  these six sets are 
     congruent to  two copies of the whole sphere,
      modulo  the countable set $D$.   
  
     Now we  perceive the above congruences,   using
     the actions $\sigma$ and $\tau$, as two sets of  colouring rules,
      with the colours $A$, $B$, and $C$.  
     We colour any point $x$ according to the colours
      of $\tau^{-1} (x)$ and $\sigma (x)$. Suppose that 
     $\tau^{-1} (x)$ is coloured  already  by $B$
     and $\sigma (x)$ is coloured already by  $C$; 
      the $\sigma$  part of the rules say that
     $x$ should be coloured
     $A$ and the $\tau$  part of the rules say that $x$ should be coloured
      $C$.  How could this conflict be reconciled?
       One way 
       would be to move further afield to determine the colour of $x$.
        For example, one could specify  a new colour $E$ that can exist at
     only one point in every orbit of $\mathbb{Z}_2 * \mathbb{Z}_3$, and
     then colour the rest of the orbit by  $A$, $B$, and $C$
     according to the group element needed to move to the point in question
     from the representative coloured  $E$ (as was done in
     the original construction of the Hausdorff paradox). After doing this,
     one could then  assign the subset $E$ to the appropriate colour
      in $A$, $B$, or $C$ as needed. 
     To insure that there is one and only one representative
    coloured  $E$ in every orbit  we could 
    require that $gx$ in $E$ for any $g\in \mathbb{Z}_2 * \mathbb{Z}_3\backslash \{ e\}$ implies that $x$ cannot be in  $E$ and
    if $gx$ is not in  $E$ for any $g\in \mathbb{Z}_2 * \mathbb{Z}_3\backslash \{ e\}$ then $x$ must be in  $E$.   However for our purposes
     such a set of   rules based on the colour $E$ 
    is  not satisfactory, 
    not only because an extra  colour has been
     added.  More  importantly such a colouring rule   would be 
     infinite in character.
     The determination of  membership in $A,B,C,$ or $E$  would involve 
     an unbounded number of group elements.

     We are interested only in  colouring rules that  are {\em finitary}
      in character, meaning that there are   finitely many colours and
      the  colour of a point
      is determined by    finitely many group or
      semi-group elements. For the rest of this paper
      we will assume that all colouring rules are finitary. 
      Of special interest are rules that always allow for
       some  assigned colour (unlike the above  Hausdorff rules); we say that such  
        colouring  rules have   {\em rank one}.
       Rank one means that 
       no matter how other point  are coloured,
       there is always an acceptable colouring for  $x$.
           Rank one rules  relate  measure theoretic
           paradoxes to problems of optimisation, especially
           when that optimisation is limited in scope. A determination
    of one best colour (or a subset of equivalently optimal colours)
    according to finite many  parameters is
    within the grasp of an automaton.

    Furthermore we  allow for colouring  rules that  depend on the 
     position in the space $X$. For example, let  $X$ be $ \{ 0,1\} ^G$
    where the group or semi-group  $G$ acts on $X$ by the
    canonical right action:  $(g(x))^h = x^{hg}$.
    There could be two distinct set of rules,
    depending on whether the coordinate  $x^e$ is equal to $0$ or $1$. 
      We  identify a  special class  of colouring 
      rules: those  that do not depend on  position,
      called \emph{stationary rules}, with
      the others called  \emph{non-stationary rules}. The distinction of
      stationary vs non-stationary
      plays an important role in this paper,
      as there are  spaces for which it is relatively easy to
      find  a paradoxical
      rank one  non-stationary  rule,  yet  it
      is unknown whether there are rank one stationary  paradoxical 
      rules with the same number of colours. As we prefer
      colouring rules that can be preformed by an automaton,
      it should be easy to determine
       which part of the rule   should apply.  
       In general, we requires that the change in rules should
      be measurable with respect to the original probability distribution.
      But in all example presented here 
       the spaces are Cantor sets and the
    domains for  the constant parts of the 
       rules are   cylinder sets. 
      
         Ramsey Theory is a weakly related topic. A  non-amenable group $G$
          could
          act on
          two probability spaces $X$ and $Y$ such that every element in $G$ is
           measure preserving in both spaces.  
          Assume there is  a paradoxical colouring rule for the space
          $X$ using the group $G$ and  $Y$ possesses   a finitely
          additive $G$-invariant extension  defined on all subsets.  
           The consequence would be that these  rules
           cannot be realised in any way on the space $Y$, though they
           would not be  contradictory in a logical  way that is independent of
            the space.
            The failure of the colouring rule on part of such a  space $Y$  is
            analogous to the existence of a forbidden structure
            with  Ramsey Theory.  The main difference with Ramsey Theory is that
            the relevant
            relations between the vertices to be
             coloured (e.g. with edges or arithmetic progressions) are
             not locally bounded, while  our colouring
             rules  are so bounded. \vskip.2cm

 The above mentioned reversal of orientation, the primacy of the colouring rules, is related to other approaches.

     A proper colouring is a colouring of the edges of a graph such that no two vertices connected by an edge are given 
     the same colour. Placed into our context of a group action,
     the edges of the graph are defined by the generators of a group acting on 
     the probability space, hence in general orbits isomorphic to
      the  Cayley graph.
     The primary concern  with Borel colouring
     (the questions raised in [4]); is what is 
     the {\em Borel} chromatic number,
     the least number of colours needed for a proper colouring
     such that each colour class  defines a 
     Borel measurable set. Furthermore
     there is interest in the $\epsilon$-Borel chromatic number -- 
     for every positive $\epsilon$ the least number of colours 
     needed for a proper colouring when the colours define Borel
     measurable sets and one exempts a measurable set of size 
     $\epsilon$ from the proper colouring requirement.
     Notice  with finitely many generators that 
     proper colouring defines a colouring rule in our context. Given
     a free action and  a common degree to
      the Cayley graph, 
     unless the number of colours is larger than this  common degree, proper 
     colouring defines a colouring rule that is not of  rank  one.
     The colouring rule of the Hausdorff paradox is a 
     restriction of the rule requiring only a
      proper colouring with three colours. It adds only two requirements
     to the rule of proper colouring, that the sets cycle $A,B,C$ in  
     the direction of $\tau$, rather than in the direction of
     $\tau^{-1}$, and that $B$ opposite $C$ (with respect
     to $\sigma$) is not allowed.  It is clear that the 
     number of colours plays a central role in proper colouring,
     that colourings into Borel sets is made much easier when the number 
     of colours increases. Likewise
     the number of colours plays a central role with paradoxical colouring rule; both in making measurable colouring easier but also  in making it 
     easier      to transform a higher rank
     paradoxical colouring rule to one that is rank one and mimicks
      the former rule.

      With a  one dimensional compact continuum of colours (a probability
       simplex determined by two extremal  colours) 
       and a Cantor set (see conclusion more of the  connection to this 
        paper), Simon
       and Tomkowicz [5]  demonstrated 
       a  colouring rule  using a semi-group action for which
       there is no $\epsilon$-Borel colouring for sufficiently
        small $\epsilon>0$ 
       (meaning a
        colouring function  that is Borel measurable
        and obeys the rule in all but a set of measure
        $\epsilon$), though there is a colouring that
     uses almost everywhere the two extremal  colours  of the continuum.

        The rest of this paper is organised as follows. The second
        section introduces the basic definitions and
        some important initial results. The third section demonstrates
        some  examples of paradoxical colouring
        rules.
        The fourth section looks at another example using
         semi-group action. 
        The fifth  section applies 
        inclusion-exclusion as  it pertain to
        probability calculations and a paradoxical colouring
        rules using the group generated freely by two generators.
        The  concluding
        section explores open questions and directions for further study.

\section{The Basics} 
Let $X$ be a probability space with
$\cal F$ the sigma algebra on which a
probability measure $m$ is defined. Usually $X$ will have a topology and
${\cal F}$ will be 
 the induced Borel sets. We say that the
relation $R \subseteq X \times X$ is \emph{admissible} if:\\
 \indent $(i)$ for any $x \in X$ there are only finitely many elements $x_1,...,x_k$ of $X$ such that $(x,x_i) \in R$ for $i=1,...,k$.\\
 \indent $(ii)$ for almost all  $x \in X$, $R$ is irreflexive.\\

 We will call the elements $x_1,...,x_k$ appearing
 in the $(i)$ the \emph{descendants} of $x$; and we say
 that $x_i$ is in the $i$-th position of $x$. \\

 \indent Let $R\subseteq X\times X$ be an admissible relation. 
 Let $A= \{ A_1, \dots , A_n\}$ be a finite
 set of colours. We say that $F: X\times A^k \rightarrow
 2^ A \backslash \{ \emptyset\}  $
 is a colouring \emph{rule} of rank one  on $X$
 if \vskip.2cm  (1) for every choice of
 of an element $b$ in $A^k$ the  $F (x,b)$
  is a ${\cal F}$  
  measurable function in $x$, and \vskip.2cm
  (2) all descendants $x_1, \dots, x_k$
   are measure
   preserving transformations   as functions of $x$,
   with respect to ${\cal F}$ and $m$. \vskip.2cm    
   The colouring rule is {\em stationary} if
   $F(x,b)$ is determined only by $b$. The colouring rule
   is {\em continuous} if the probability space $X$ has
   a topology, ${\cal F}$ are the Borel sets, and 
     for every choice of
 of an element $b$ in $A^k$ the  $F (x,b)$
 is a continuous function of  $x$. All of our examples are continuous, with
 the space $X$ a Cantor set, and the different  applications
 of the rule determined
  by cylinder sets. 
   A colouring rule of rank $n$ is a collection of
   $n$ colouring  rules of rank one such that they cannot
    be expressed as a collection of fewer colouring rules of rank one.\\

    A colouring $c: X \rightarrow  A= \{ A_1, \dots , A_n\}$ {\em satisfies}
    a colouring rule $F$ of rank one at $x\in X$ if almost
     everywhere (with respect to $m$) it follows
     that $c(x)\subseteq F (x, c(x_1), \dots c(x_n))$.
     It satisfies a colouring rule of higher rank if it satisfies all
     the correspondences $F$ defining the  colouring rules of
     rank one that make up  that colouring rule of higher rank.\\

  \indent We say that a  rule  on $X$
  is \emph{non-contradictory} if it is satisfied by some colouring. 
  It is paradoxical if it is non-contradictory and for every 
   finitely additive extension $\mu$ of $m$ for which the descendants are still
   measure preserving there is no  colouring satisfying the rule
   almost everywhere such that the sets defined by the colours are measurable with
    respect 
    to $\mu$. \\

 \section{First Examples} We need to  establish  that
 there are paradoxical  rules of rank one. The first example
 is that of a rank one non-stationary rule with five colours on a space
 that may fail to have any paradoxical stationary rule of rank one using 
  the same number of  colours.
  The second and third  examples are
  that of a stationary  rank one  paradoxical  rules which
  mimick the Hausdorff rules and paradox. With the second example
  this is done by doubling the space and with the third example by
  squaring the number of colours. We are not so interested in
  the second and third examples because they represent a kind of
   cutting of 
   a Gordian knot. As stated in the conclusion,
   we would like to find some way to differentiate between these two
    examples and the others.\\

 {\bf Example 1:} Let $G$ be $\mathbb{Z}_2 * \mathbb{Z}_5$, with
 $\mathbb{Z}_2$ generated by $\sigma$ and $\mathbb{Z}_5$ generated
 by $\tau$. The space $X$ is $\{ 0,1\} ^G$ and 
 there are five colours, $A_1, A_2, A_3, A_4, A_5$. All arithmetic is
 modulo 5 and we present the rules with some simplifying abuse of notation.
 The rule is deterministic, meaning it is a function. We
  assume that $\tau^{-1} (x)$ is coloured  $ A_{i}$.
  If  either $x^e=1$, 
  $\tau ^{-1} (x)$ is coloured $ A_5$  or $\sigma (x)$ is coloured $ A_1$ then
  $x$ is coloured  $ A_{i+1}$. Otherwise, if none of
  these three conditions apply,  $x$ is coloured $ A_i$.\\

  {\bf Theorem 1:} The colouring rule of Example 1 is paradoxical. \\

  {\bf Proof:} 
  If the colouring  is well defined and follows the colouring rule,
  there are only two possibilities: 
  every  $\tau$ cycle $(x, \tau (x), \tau^2(x), \tau^3(x), \tau ^4(x))$ is
  all of the same colour (using  a single colour $A_i$ for some $i$) or
  it contains all five colours. This allows us to make two independent
  calculations concerning the proportion of the set coloured  $A_1$
  (a proportion 
  we assume exists under the assumption of finitely additive  measurability 
   and $G$ invariance).

  First we estimate from below  and outside 
   the $\tau$ cycles the probability of the set coloured $A_1$. In every
  $\tau$ cycle, there are a certain number of points whose $e$ coordinates
  are $0$. If just one of the $e$ coordinates is $1$, it
  is necessary for the cycle to include all five colours.
  For each possibility for the number of points whose $e$ coordinate is $0$,
  the probabilities of these cycles
   are based on $\frac 1 {2^5} = \frac 1 {32}$ and  the
  binomial expansion. We consider 
  the points opposite to  the  cycle, namely
  the points $\{ \sigma (x), \sigma \circ \tau (x), \sigma \circ
  \tau^2(x), \sigma \circ \tau^3(x), \sigma \circ \tau ^4(x)\}$, and see
  how many of them must be coloured  $A_1$.  These calculations are based
   on $\frac 1 {160}$, with  $160= 32 \cdot 5$. 
   The probability of  either one or no such  points in the cycle
    with $0$ as its
    $e$-coordinate 
  is $\frac {1+5}{32}$, and in that case there need not be any points
  opposite the cycle in $A_1$ (as the one place $x$ in
  the cycle where $x^e=0$ could be the one place 
    that is coloured  $A_5$).
  Following the same logic,
  if there are exactly $k$  points in the cycle whose $e$ coordinate is $0$,
   with $2\leq k \leq 4$, then
   there must be at least $k-1$ points opposite the cycle
   coloured by   $A_1$.  If $k=5$ then
   of course the entire cycle may be coloured with just one
   colour, not requiring any opposite point to be in $A_1$.
    The proportion of the space taken up by the  set $A_1$ must be at least
    $\frac 15 ({5\choose 2} + 2 \cdot {5\choose 3} + 3 \cdot {5\choose 4} ) =
    \frac {10} {160} + \frac {20}{160} + \frac {15} {160} = \frac {9} {32}$.

    Second  we estimate from above and inside a $\tau$ cycle
     the probability of the set coloured $A_1$.
    The probability of a $\tau$ cycle having the same colour throughout cannot
    exceed $\frac 1 {32}$, the probability that all members of the cycle
    have $0$ for the $e$ coordinate. Otherwise at most one-fifth of the
    cycle can be in $A_1$. This implies that the proportion of the set
     coloured $A_1$ in
    the space cannot exceed $\frac 1 {32} + \frac 1 5 \cdot \frac {31}{32}=
    \frac 9 {40}$. Therefore the proportion, if it exists, must be below
     $\frac 9 {40}$ and above $\frac 9 {32}$, a contradiction. 

    A  non-measurable  colouring  satisfying  the above rule is  
     essentially the 
    same as that for the Hausdorff paradox:
    all  $\tau$  cycles
    go through all five colours and all but one of the five points opposite
     to any cycle 
     opposite are coloured  $A_1$ (and the $e$ coordinates are ignored). 
    \hfill $\Box$ 
\vskip.2cm 

    {\bf Example 2:} Let $G= {\mathbb{Z}_2 * \mathbb{Z}_3}$  and
    $X= \{ R,S\}  \times \{ 0,1\} ^G$, with half probability
    given to both  $\{R\} \times \{ 0,1\}$ and  $\{S\} \times \{ 0,1\}$.
     We define $\rho$ to be the measure preserving  involution that switches between $(R,x)$ and $(S,x)$ for every $x\in \{ 0,1\}^G$. 
    As before, 
    $\mathbb{Z}_2$ is generated by $\sigma$ and  $ \mathbb{Z}_3$ is generated
    by $\tau$. We use the three colours $A, B,C$,  the same from 
     the Hausdorff rule. 
    The new rule is  as follows. It   assigns  to
    every  $(S,y)$ the same colour as that of $(R,y)$.
     The rules for colouring the points of the form
    $(R,y)$ are based  on the Hausdorff 
      rules  and  
     are very easy to describe.  We have described
    above the rank two rules defining  the Hausdorff paradox,
    determining the colour of $y$ according
     to the colours of $\sigma (y)$ and $\tau^{-1} (y)$. If the colour
     of $(S,y)$ follows the Hausdorff rule with respect
      to $(R,\sigma y)$ and $(R, \tau^{-1} y)$ (after dropping the $R$ and $S$),  then
    colour $(R,y)$ the same colour as that of $(S,y)$.
    Otherwise colour $(R,y)$ differently than
     $(S,y)$, with either other way valid.    
    In any  colouring satisfying the rule the points   
    $(R, y)$ and $(S,y)$ are given the same colour, implying
    that the Hausdorff rules are being  followed on both
    copies of $\{ 0,1\}^G$.\\

    {\bf Example 3:} 
    Return to the uncopied space $ X= \{ 0,1\} ^{\mathbb{Z}_2 * \mathbb{Z}_3}$. 
    By using nine colours, namely the pairs $(A,A)$, $(A,B)$, $\dots$,
    we can create an analogous rule to that of Example 2,
    where the second colour of  $x$ is the 
    copy of first colour of $\sigma x$, and the first colour of $x$ agrees
    with the second colour of $\sigma x$ if and only if the
    second  colour of $\sigma x$ is the correct colour for $x$ when
     following the Hausdorff rule with respect to the first colours of $\sigma x$ and $\tau^{-1} (x)$. This shows that
    there is a stationary rank one paradoxical colouring rule that recreates
    the Hausdorff paradox after dropping the second colour.

       \section{A colouring  rule that is paradoxical using a   semi-group}
               {\bf Example 4:}
               Let $G = \mathbb{N}_0 * \mathbb{Z}_2$ be the free product with generators $T$ and $\sigma$,
       respectively, where $\mathbb{N}_0= \{0,1,2,\dots \} $, meaning that
       only non-negative powers of $T$ are used. There are two colours, $A$ and $B$. 
       Let   $X = \{0,1\}^{G}$ be the Cantor space with the canonical measure
       $m$. Consider
       the following colouring  rule $Q$ involving two colours $A$ and $B$; ({\em opposite to}  means using the other
        colour)
        :\\

 $(i)$ if $T x$ is coloured $A$ then colour $x$   opposite to the colour of  $\sigma x$;\\
 \indent $(ii)$ if $T x$ is coloured $B$ and $x^e = 0$ then colour $x$ the same as  colour of   $\sigma x$;\\ 
  \indent $(iii)$ if $T x$ is coloured $B$ and $x^e = 1$ then colour $x$  opposite to  the colour of  $\sigma x$;\\ 

  It seems  that this colouring rule $Q$
  must be more than paradoxical, namely contradictory. Every pair $x,y\in X$ determine four pairs of  points $z$ and $\sigma x$  such that
  $Tz =x$, $T \sigma z =y$, namely the four possibilities of $z^e= 0,1$ combined with $(\sigma z)^e=0,1$. Furthermore,
  from the measure preserving properties of $T$ and $\sigma$ the distribution $m$  on $X$ is determined  independently by
   those four choices of pairs and a choice 
   of $x$ and $y$ in $X \times X$ in the $T$ and $T\sigma$ positions,  according to the product measure $m\times m$.  If both $x$ and $y$ are coloured $B$, then half of the four
  pairs of  possibilities of $z$ and $\sigma z$ cannot be coloured at all, since if $z^e=0$ and $(\sigma z)^e=1$ (or vice versa)
  there is no way to colour one or the other of these two points. 
  Also,  if $x$ is 
  coloured $B$ and $y$ is coloured $A$ 
     then either the    $z$ satisfying  $z^e=0$ or its twin   $\sigma z$  cannot be coloured. 
   Hence from the measure preserving property of $T$ and $\sigma$, any positive measure $r$ given to
   the subset coloured $B$ results in at least $\frac r4$ of the space being uncolourable according to the rule. On the other hand, from the
   probability of at least $\frac 12$ that $(\sigma x)^e$ is equal to $1$,
   the subset coloured $A$ cannot exceed a size of $\frac 56$, meaning that the colouring rule cannot hold in  at least $\frac 1 {24}$ of the space. Clearly a  measurable
   colouring is not possible by any interpretation of
   measurability. But even without assuming that
    the sets coloured $A$ and $B$ are measurable in any way, it seems that  if they are to be coloured, then
     it  must  be  in a  measurable way. 
     If the colouring rule holds  almost everywhere relative to the Borel measure (with or without measurability of any kind to the colouring),
     from the above logic the whole set $B$ presents a problem to  the colouring of $T^{-1} (B)$--  therefore $B$ should be contained in a set
   of Borel  measure zero. But this means that $B$ is  measurable and of measure zero with respect to the completion of the Borel measure, a contradiction.

    Indeed showing that
    this colouring rule $Q$ is  paradoxical (e.g. also  satisfiable)
     is different from the previous examples. Rather than first showing that the rule cannot be satisfied in any measurable way and then showing
   satisfaction of the rule by some colouring, we reverse the process and show first the satisfaction. In order to do this, we must use only part of the set $X$ and redefine the Borel  measure
   to apply only to this subset. In this special  subset of $X$, the set coloured $B$ will never be in the image of $T$,
   hence the problem presented above will be avoided.

We define the functions $\alpha : X \rightarrow X$ and $\beta : X \rightarrow X$ by
$\alpha (x) ^e + x^e=1$ (modulo 2), $\beta (x) ^ {\sigma} + x^{\sigma} = 1$ (modulo 2), and
otherwise $\alpha (x)^g = x^g$ for all $g\not= e$ and $\beta (x)^g = x^g$ for all $g\not= \sigma$. In other words,
$\alpha$ switches the $e$ coordinate and $\beta$ switches the $\sigma$ coordinate. Both $\alpha$ and $\beta$ are measure preserving transformations.
We include $\alpha (x)$ and $\beta (x)$ among the descendents of $x$, although they will play no rule in the colouring of $x$.
We extend the semi group $G$ to $\overline G$, the semigroup generated by $G$ and $\alpha$ and $\beta$.  From the above analysis,
 for every pair of  subsets $D,E\subseteq X$, $\alpha$ and $\beta$ map the set $T^{-1}(D) \cap \sigma T^{-1} (E)$ back to itself.

 We  use the same  approach that has been applied to Wiener processes and the restriction to the paths that are continuous. With this approach, $\Omega$  is a topological  space 
 with a regular Borel probability measure  and  
 $C$ is a  non-measurable subset of outer measure one. The subset
  $C$  inherits the original structure of the  Borel measure through intersection, so that  any Borel set $B$ of $\Omega$  defines a measurable subset
   $B\cap C$ of $C$ whose measure $m (B\cap C)$ is the same as $m (B)$ (see [2] for more details).  

We will show first that there is a  colouring following the rule $Q$ on a subset $X'$ of outer measure 1. Except for the outer measure, the  method  
is similar to that used in Simon and Tomkowicz [5]. The measure
    preserving property of the $T$,  $\sigma$, $\alpha$ and $\beta$ 
    remain after such a restriction  to the set $X'$.

We say that a set $S(x) \subseteq X$ is \emph{the semigroup orbit of x} if it has the form\\
\indent $\{y \in X: y = U(x) \ \textrm{and} \ U \ \textrm{is a semigroup word in} \ T, \sigma \}$. A subset $A$ of $X$ is {\em pyramidic} if for every $x\in A$
the semigroup orbit of $x$ is also in $A$. The link $R(y)$  of a point $y$ is
$\{ U^{-1} (y)\ | \ U \mbox { is an element of G} \}$. The semi-group orbit
   of a point is a countable set, while its link is an uncountable set.  However due to the measure preserving property of all
   elements of $G$ and that there are countably many members of $G$, all links are sets of measure zero (following from
    the atomless character of the probability space).\\

 \indent We define the \emph{outer  measure} $m^{*}$ of $Y \subseteq X$ by 
 $$\inf \{m(L): Y \subseteq L: L \ \textrm{is an open set of} \ X\}.$$
 
 Assuming the continuum hypothesis or weaker Martin's axiom we have the following theorem:\\

 {\bf Theorem 2:} There exists a pyramidic  subset $X'$ of outer measure one
 and a 
 colouring $c$ of $X'$  that satisfies the rule $Q$
 such that $T(X') \subseteq A$  and
  the set of points coloured $A$ has outer measure $1$. \\

 {\bf Proof:} We will be proceeding by transfinite induction. 
 First we order all uncountable compact subsets of $\{0,1\}^{G}$ with positive
 measure $m$ into a transfinite sequence $\{K_{\alpha}\}$. Since $\{0,1\}^{G}$
is a separable space and the open sets have
a countable base,
there are continuum many compact sets of $\{0,1\}^{G}$.
Pick now a point $x_0 \in K_0$ and colour it $A$.
Then consider the semigroup orbit $S_0$ of $x_0$ and colour it
in the way that $TV(x_0)$
is always coloured $A$, where $V$ is a semigroup
word in $\sigma$ and $T$. Clearly, the  colouring
of $S_0$ satisfies the rule $Q$.\\
\indent Suppose now that we have coloured the semigroup
orbits $S_0,...,S_{\beta}$ for some ordinal $\beta < \mathfrak{c}$.
 Clearly, the set $T_{\beta} = \bigcup_{i  \leq \beta}S_i$ has
 the cardinality $\mathfrak{m}$ less then the continuum $\mathfrak{c}$ 
 (this follows from the fact that any cardinal number is the
 least ordinal number with the given cardinality).
 \\
 \indent Now, given a point $y \in T_{\beta}$,
 we consider the set of links $R(y)$ and then the set
 $W_{\beta} = \bigcup_{y \in T_{\beta}} R(y)$.  
 As the links are sets of measure zero,
 the set  $W_{\beta} = \bigcup_{y \in T_{\beta}} R(y)$ is 
  also a set of measure zero.
  Now take a point $x_{\beta+1}$ in $K_{\beta+1} \setminus W_{\beta}$.
  Clearly, for any word $U \in G$
  we have $U(x_{\beta+1}) \notin T_{\beta}$.
  Hence we colour $x_{\beta +1}$ in $A$ and we colour the remaining points
  in the semigroup of $x_{\beta +1}$ in the way to satisfy the rule $Q$.      
  Therefore for any compact set $K_{\alpha}$ there
  is a point $x_{\alpha} \in K_{\alpha}$ coloured $A$.
  If $A$ failed to have outer measure one,
  there would have to be a cover of $A$ by
  open sets whose measures add up to strictly less than $1$,
  hence there would have to
  be a compact set of measure strictly greater than $0$ with an empty intersection with $A$, something
   we have excluded through the construction. 
  Finally we define  $X'$ as the union of the sets coloured $A$ and $B$.  $\Box$\\ 

  Still there  does seem to be   something wrong with the
  above theorem.
  With $T$ measure preserving, one could
  think that  $X'= T^{-1} T (X')$
  should be the same measure as that of $T(X')$,
  which is the subset coloured $A$. This would
 suggest, after an abuse of notation,
 that $B=\sigma (A)$ is also of measure one, a contradiction.
 But we never  imply in the context of intersection with  $X'$
  that $T(X')$ is a measurable set, and indeed  the measurability
  of $X'$ does not imply the measurability of
  $T(X')$ (as it is an uncountable to one function).\\

  {\bf Corollary:} The colouring rule $Q$ applied to $X'$ is paradoxical,
  given that $\alpha$ and $\beta$ are descendants and 
   the finitely additive
  extensions are  $\overline G$ invariant.\\

   {\bf Proof:} 
   Having defined  $X'$ as above, we have a colouring into $A$ and $B$
   that satisfies the colouring rule. Now remove the initial  partition of $X'$
   into the two colours $A$ and $B$ and
   consider any  colouring according to the rule $Q$
   (not necessarily in the same way as we used
   to construct $X'$). We assume that this new colouring
   is measurable with respect to a finitely additive
   measure $\mu$ which is $\overline G$ invariant and according to
   this new colouring we call $\overline A$
   the subset coloured $A$ and $\overline B$
   the subset coloured $B$.   With  both $\alpha$ and $\beta$ mapping 
   $T^{-1}(\overline B) \cap \sigma T^{-1}(\overline B)$
   to itself and $T^{-1}(\overline B) \cap \sigma T^{-1}(\overline A)$ to itself, from the  $\alpha$ and $\beta$ invariance
   we must conclude (for both  sets)
   that the four possibilities for the
   $e$ and $\sigma$ coordinates are equally likely (according to $\mu$). 
   Following the  same arguments as above,  
   at least $\frac 1 {24}$ of the space cannot
   not be coloured according to the rule $Q$, a
    contradiction. \hfill $\Box$\\

{\bf Note:}    An alternative approach would be to construct $X'$ using
     orbits and the links defined  according to larger group 
     $\overline G$. In the proof of Theorem 2, one could  try to make
      the set coloured $A$  equal to the image $T(X')$. 
 
 \section{Inclusion-Exclusion principle}

     In this section
     we will present an application of the Inclusion-Exclusion principle 
     to paradoxical rules.  Let $\mu$ be any finitely additive measure. 
     The inclusion-exclusion rule  for finitely
     many sets $A_1,...,A_k$ can be phrased as:
 $$ \mu(\bigcup_{i=1}^{k} A_i ) = \sum_{j=1}^{k}(-1)^{j+1} \Bigg{(} \sum_{1 \leq i_1 \leq ... i_j \leq k} \mu (A_{i_1} \cap ... \cap A_{i_j}) \Bigg{)}.$$

Consider now some  sets $A_1, A_2$ and $A_3$ contained in $X$
and such that $\mu(A_1) + \mu(A_2) + \mu (A_3) = 1$.
Then, by the Inclusion-Exclusion principle and the
fact that $\mu(X) = 1$, we get that\\

$(*)$ \  $\mu(X) - \mu(A_1 \cup A_2 \cup A_3) = \mu(A_1 \cap A_2) +
\mu(A_2 \cap A_3) +\\
 \mu(A_1 \cap A_3) - \mu(A_1 \cap A_2 \cap A_3) = 0$.\\

 Clearly, the same equality works, mutatis mutandis, for any $k$ subsets $A_1,...,A_k$ of $X$ such that $\mu(A_1) +... + \mu(A_k) = 1$.
 Of particular interest is when  $B_1, \dots , B_k$ partition a probability space and $A_i = \sigma_i B_i$ for some
 measure preserving transformations $\sigma_1, \dots , \sigma_k$.
 If  we can show that
 $\mu(A_i \cap A_j)=0$ for all $i\not=j$ then any argument that
  $X\backslash (A_1 \cup \cdots \cup A_k)$ should
 have positive measure would be a contradiction to the measurability assumption.\\

The following   example demonstrates a  
 non-stationary paradoxical rule 
 with four colours for the space $\{ 0,1\} ^{\mathbb{F}_2}$  using
 the inclusion-exclusion principle. As $\mathbb{F}_2$ is isomorphic to a subgroup of $\mathbb{Z}_2*\mathbb{Z}_3$, this example shows also that
 there is a paradoxical
 colouring rule for 
 $\{ 0,1\} ^{\mathbb{Z}_2 * \mathbb{Z}_3}$  with four colours.
 With $\sigma$ the generator of
 $\mathbb{Z}_2$ and $\tau$ the generator of
 $  \mathbb{Z}_3$, this follows from the fact that 
 $ \tau \sigma \tau^2\sigma$ and
  $ \tau^2 \sigma \tau\sigma$ generate of a subgroup of $\mathbb{Z}_2 * \mathbb{Z}_3$ isomorphic to $\mathbb{F}_2$.

{\bf Example 5:}
 With $X= \{ a,b\} ^{\mathbb{F}_2}$, initially there will be two colours,
 $P$ and $N$ for positive and negative.
 Before  specifying the rule determining
 the colouring, let $P_i$ be the subset
 coloured positive  such that  the $e$ coordinate is $i=a,b$, and
 define $N_i$ in the same way.
 Let $T_a$ and $T_b$ be the two generators of $\mathbb{F}_2$.
 If we can create a rank one rule such  that
 the four sets $T_a(N_a)$, $T_a^{-1} (P_a)$, $T_b(N_b)$, $T_b^{-1} (P_b)$
 are  mutually disjoint almost everywhere, then automatically
 the rule must be paradoxical. This follows from
 the fact that the set of $x$ such that
 $T_a(x)^e=b , \ T_a^{-1}(x)^e=b, \ T_b(x)^e= a,\  T_b^{-1}(x)^e=a $
 is $\frac 1 {16}$ of the whole space, yet by the inclusion-exclusion
 principle the measure of this set must be $0$.
 
 To  define the  rule, we will create four colours by
 splitting the initial
 two  colours $P$ and $N$ into
 four colours, $P^u$, $P^c$, $N^u$ and $N^c$, the $C$ or $U$ for
  whether the vertex is ``crowded'' or ``uncrowded''.
  Given any particular colouring, we place an arrow
  from every $x$ to one of its four neighbours. If $x^e=a$, we place an arrow
  from $x$ to $T_a (x)$ if $x$ is coloured either $P^u$ or $P^c$
  and otherwise  place an arrow
 from $x$ to $T_a^{-1} (x)$.
 Likewise if $x^e=b$, we place an arrow from $x$ to $T_b(x)$
 if $x$ is coloured $P^C$ or $P^U$ and otherwise place an arrow
  from $x$ to $T_b^{-1} (x)$. 
  Now in defining the  rule, we have to decide for every $x\in X$
  whether  to colour it $P$  or $N$
   and whether to colour it crowded $C$  or uncrowded $U$.
  The rule is relatively simple. If $x^e=d$ with $ d \in  \{ a,b\}$  and  $T_d(x)$
  is uncrowded and $T_d^{-1} (x)$ is crowded then colour $x$ with $P$.
  Likewise of $x^e=d$  and  $T_d(x)$ is crowded and $T_d^{-1} (x)$ is
  uncrowded, colour $x$ with $N$.
  Otherwise if  the appropriate adjacent vertices are both
  uncrowded or crowded then colour $x$ however one wants.
   A vertex is coloured ``crowded'' or  $C$ 
  if there are two or more arrows pointed toward it and it is coloured
   ``uncrowded'' or  $U$ if there is one or no arrows pointed toward it.\\

   {\bf Theorem 3:} The colouring rule of Example 5 is paradoxical.\\

   {\bf Proof:} 
   We define the degree of a vertex
   as the number of potential arrows that could be pointed
   toward this vertex. In other words,
   if $(T_a (x))^e=a$, $(T^{-1} _a (x))^e=b$, $(T_b (x))^e=b$ and $(T_b^{-1} (x))^e=b$ then the degree of $x$ is $3$.
    If  $(T_a (x))^e=b$, $(T^{-1} _a (x))^e=b$, $(T_b (x))^e=a$ and $(T_b^{-1} (x))^e=a$ then the degree of $x$ is $0$. 
    As mentioned before,
    the degree zero  vertices take up $\frac 1 {16}$ of the space.  
    Our claim is that with a colouring satisfying the rule
    almost all vertices are uncrowded.
    This implies that the colouring cannot be measurable with respect
    to any finite extension for which the group $\mathbb{F}_2$ is measure
     preserving, 
   because  at least $\frac 1 {16}$ of the vertices  cannot have
   any arrows points toward them. By one accounting
   there is an arrow exiting every vertex but by another
   accounting the average number of arrows coming in to vertices
    must be no more than
   $\frac {15} {16}$.

   Now let us assume that $x$ is a crowded vertex and see what is necessary
   to maintain this situation in a colouring satisfying the rule.
   There must be two distinct vertices $y$ and $z$ such that
   there is an arrow from $y$ to $x$ and an arrow from $z$ to $x$. Lets
   focus on just one of them, without loss of generality the $y$. As the
   colouring rule is satisfied,
   the existence of an
    arrow from $y$ to $x$ implies that the vertex
   opposite to  $x$ from $y$   is another crowded vertex, call it $w$.
   As the arrow
   is already defined from $y$ to $x$
   it means that there are at least  two arrows
   pointed inward to  $w$ that do not start at $y$. Letting $v_1$ and $v_2$
   be two of those vertices, let $u_1$ and
   $u_2$ be the vertices for which there could have  be an arrow
   from $v_i$ to $u_i$, however instead the arrow was
   from the $v_i$ to $w$. We recognise by induction
   the existence of a chain of backwardly directed
    arrows, starting at $w$, moving to
    $u_1$ and $u_2$ and beyond, such that the induced graph is binary,
    has two branches at every stage. (If there are three
    such branches in some places, we could reduce
    to the existence of a binary tree). Each of these
   vertices of the binary tree  has degree at least three. 
   Now let $p$ be the probability of there existing
   such an infinite chain, the probability relative
   to the start at a   vertex like $y$ moving in the direction  away
    from  $x$. 
As the space is
defined homogeneously (that the probabilities for $a$ or $b$ are
independent regardless of  shift distances)
we can calculate $p$ recursively.  
There are two possibilities, the next vertex $w$  could be of degree
three and the chain of backward arrows
continues with these two adjacent vertices on the other side of
$y$,  or $w$ is degree four and 
the chain continues with at least two of these three vertices on
the other side of $y$.
In the first case the conditional  probability that $w$ is of degree three
is $\frac 38$ (conditioned on the move from 
      $y$ to $w$) and then a  continuation of the chain indefinitely 
happens with probability $p^2$. In the second case,
the conditional  probability of degree four is
$\frac 18$ and the probability of continuation $3p^2 - 2p^3$
(three choices for  the two next vertices  minus the
possibility, counted twice,
that continuation is possible in all three directions).
We have the
formula $p = \frac 38 p^2 + \frac 18  (3p^2 - 2p^3)= \frac {3p^2 - p^3} 4$.
Factoring out the $p=0$ solutions, we are left
with $p^2 -3p +4$, which has no real solutions.
We conclude that the  stochastic structure of degrees  does not allow
     for crowded points to exist in more than a set of measure zero.

     Now we show (using AC) that there does exist a  colouring satisfying
      the rule. For every orbit choose a representative $x$ and label every other vertex in this same orbit as $gx$ by the group element $g$ used
      to travel from $x$ to $gx$. Every group element $g$ has a length, the
      minimal  number of uses of $T_a$, $T_a^{-1}$, $T_b$ and $T_b^{-1}$
      used to construct $g$ (where the length of
      the identity is zero). Let the length of a vertex $gx$ in the orbit be the length of $g$.   Colour $x$ with
      either colour. Colour all  vertices of length $1$ next, then all
      vertices of length $2$, and so on. At every stage of the process, a vertex of length $l$ is adjacent to one vertex
     of length $l-1$ and three of length $l+1$. Therefore from any vertex $y$ of length $l$ one can always point the arrow toward a vertex of length $l+1$, regardless of the value
       of $y^e$. There can be no other arrow toward $y$, as the three other points adjacent to $y$ are all  of length $l+2$. \hfill $\Box$

 \section{Further study}

 {\bf Question 1:} Given  any probability space $X$ and a 
 finitely generated measure preserving group $G$ acting on $X$,
 let $P$ be the colouring rule
 that   requires only
 that the colouring must be proper. Is there such an $X$ and $G$
  such that the rule $P$ is  paradoxical? 
  \\
  
    A  colouring rule with $k$ colours
     is (stationarily) \emph{essential} if
    it is paradoxical of rank at least two and there exists no (stationarily)
     paradoxical rank one finitary rule with $k$ colours whose satisfaction
     implies the satisfaction of the original rule.
     A space $X$ is (stationarily) essential
      with $k$ colours
      if there is no (stationarily) paradoxical rank one colouring  rule
       with $k$ colours defined
      on $X$, and
      yet there is  a paradoxical finitary rule with $k$ colours
       defined on $X$  of higher rank. 
 Intriguing is the difference between what can be accomplished by
 stationary and non-stationary colouring rules.\\ 

 {\bf Question 2:} Is 
  there is a non-stationary rank one paradoxical colouring rule
  on $X=\{ 0,1\} ^{{\mathbf Z}_2 * {\mathbf Z}_3}$ using
  the canonical measure preserving
  group action of
  ${\mathbf Z}_2 * {\mathbf Z}_3$ and three colours, but  no such stationary
  paradoxical rank one paradoxical colouring rule with three colours? \\
   We conjecture
   that the answer to Question 2 is yes. 
  In other words, we conjecture that this space with this group action
  is stationarily essential but not essential. If indeed it is not
  stationarily essential, there remains the question of whether
  the conventional Hausdorff rule of rank two is essential, utilising either 
  a stationary or non-stationary rule involving a larger number
   of descendants.\\ 

{\bf Question 3:}    Does there exist a space and a semi-group action
    that is essential with respect to any number
    of colours? \\

    Whenever at least one of the descendents is invertible
    and its inverse is measure preserving (the latter following from the
    former whenever the space is compact and the  measure is Borel)
     then a colouring rule
    of any rank with $k$ colours can be mimicked by a rank one colouring
     rule with $k^2$ colours in the same way as with Example 3. 
    But if none of the
    descendants are invertible, it is plausible that some paradoxical
    colouring rules are  not multi-functional in character no matter how
     many colours are allowed. \\

    There is something satisfying about Example 5 and unsatisfying
    about Examples 2 and 3 (with
       Example 1 somewhere in between).  Example  5 employs a
    stochastic process that seems to push colourings toward the satisfaction of
     the rule.
     On the other hand, the existence of  colourings witnessing the rules
     of  Examples 2 and 3
      seems to be either accidental or  contrived.
     One could perceive
  satisfaction of a colouring rule to be a kind of fixed colouring, with
  the colouring rule defining some kind of iterative process that does or
   doesn't bring the colourings closer to satisfaction of the rule. 
   Of course if any colouring rule forced almost everwhere
   (with respect to the original measure $m$) an eventually stable 
    colour in finitely many colouring stages  (with respect
    to  some initial measurable  colouring not satisfying
     the rule) then there would be measurable
    colouring solutions and the rule could not be paradoxical.\\ 

{\bf Question 4:}     Is there 
a  refinement  to  the  definition of a  rank one paradoxical rule  
  that
    identifies a credible   force toward its satisfaction?

 Suppose one had a colouring rule for $S ^G$ where $S$ is a finite set
  and $G$ is a group, and $C$ are the finite set of colours. 
  Another  structure to consider is $\{ S \times C\} ^G$,
  where we assume a random colouring start to $S^G$ inherited from
  the $C$ coordinate. We could consider how the colouring
   rule generates iterations of colourings on $\{ S \times C\} ^G$. 
 
       The \emph {deficiency} of a  colouring
        rule on $X$ is the infinum of all
        probabilities $\rho$ such that there is a
         ${\cal F}$ measurable  subset  
      of size  $1-\rho$ such that
      the finitary rule {\bf can}  be satisfied on this subset.
      The paradoxical  {\emph size} of a
      paradoxical   colouring rule on $X$
      is its deficiency. Likewise we define the paradoxical 
      size of  a colouring rule as the difference between
       its  deficiency and the supremum
      on all  $r$ such that there is a satisfaction of the rules on a
       ${\cal F}$ measurable subset of size $1-r$. 

      Do   paradoxical colouring rules of arbitrarily  small
       paradoxical sizes get   converted in an organic way 
         to full paradoxical decompositions? 
         With the Hausdorff paradox the ambiguity concerning the size of
          the set coloured $A$, that it should 
        be both $\frac 12$ and  $\frac 13$, leads to a statement of
        ``$1=2$'' for the whole space.  
        More generally
        Tarski proved that the existence of a finitely additive
        measure defined on all subsets is equivalent to the existence
        of a paradox in the conventional sense, that there are disjoint
         sets $B_1, \dots, B_l$  
         of a subset $E$ of positive measure and measure preserving
         transformations $\sigma_1, \dots, \sigma _l$ with
         $\sigma_1 (B_1), \dots, \sigma_l (B_l)$ forming two
         copies of $E$. If there is a paradoxical colouring
          rule on $X$
         using a group $G$ of measure preserving transformations, we know
         automatically that $X$ is non-amenable in character and
          thus has such a paradoxical decomposition without
          any application
          of the paradoxical rule.\\

          {\bf Question 5:} Given a colouring $c$  
          that satisfies a paradoxical colouring rule, are there 
          disjoint subsets $B_1, \dots , B_l$ of a set $E$ of positive
          measure  creating  two copies of $E$ using the
           measure preserving descendants 
           such that the  $B_i$ sets belonging to the
         smallest sigma algebra containing the original ${\cal F}$,
         the sets defined by the inverse images of the colouring $c$,
          and the shifts of these sets as so  defined by the descendents? 
\\
          
          In the definition of a colouring rule, we use that the
          descendants are measure preserving. This was a natural way
           to connect colouring rules to measure theoretic paradoxes. 
\\
           
{\bf Question 6:} 
Is there a more general definition for paradoxical
colouring rules (implying that the sets defined by the  colours cannot not
  be measurable by any invariant finitely additive measure) 
           that  weakens the measure
           preserving property of the descendants?
           \\
           
           {\bf Question 7:}
           When   the space $X$ to be coloured has a topology, is
            there any significant distinction between
             the continuous  
             and the more general  non-stationary colouring rules? 
           \\

Usually 
optimisation concerns a choice of a point in a finite dimensional
convex set.
We  generalise to a  definition of a
{\em probabilistic}   rank one   colouring
rule in the following way.
We keep the initial set up with $X$ a probability space and
 the finitely many descendants, but we change
 what is a colour. The set of colours is  
$A = \prod_{i=1} ^n \Delta (C_i)$, where each  $C_i$ is a finite set.
 The  restriction on the choice of a colour 
 is to a   non-empty  subset of $A$,
 a restriction determined by  $x$ and
 the descendents $p_1, p_2, \dots p_k \in A^k$ in a measurable
 way. The probabilistic  colouring rule is   {\em normal}  
 if  these non-empty subsets of $A$ are convex and
  as a correspondence defined by
  the $x$ and the $p_1, \dots, p_k$ is upper-semi-continuous.  We 
 define satisfaction of the colouring rule in the same way as before and
  understand by 
  measurability with respect to a finitely additive
  extension $\mu$ of $X$ such 
 that  the inverse images of any Borel set in $A$ are
  $\mu$ measurable. As before, a probabilistic  colouring rule is paradoxical
 if it can be  satisfied but not by any colouring function which
  is measurable  with
 respect to any finitely additive extension that maintains
 the measure preserving property of the descendants.
 \\
 
{\bf Question 8:} 
 Does there exist a normal  probabilistic rank one   paradoxical
 colouring rule? \\
 
    This would be of great importance to applications, as its affirmation
    would  suggest
 that there may be some very natural
 optimization problems that are solvable in some way, but not
  in any way for which one can make any  predictions.\\

\end{document}